\documentclass [11pt]{article}
\usepackage{amsmath}
\usepackage{amssymb}
\usepackage{amscd}
\usepackage{amsthm}
\usepackage{amsopn}
\usepackage{xspace}
\usepackage{verbatim}
\usepackage{amsmath}
\usepackage{amsfonts}

\newtheorem{lemma}{Lemma}[section]
\newtheorem{proposition}{Proposition}[section]

\newtheorem{corollary}{Corollary}[section]

\newtheorem{acknowledgment*}{Acknowledgment}

\numberwithin{equation}{section}

\newcommand{\ooM}{\overline{\cal M}}

\newcommand{\be}{\begin{equation}}
\newcommand{\ee}{\end{equation}}
\newcommand{\bd}{\begin{displaymath}}
\newcommand{\ed}{\end{displaymath}}

\newcommand{\eps}{\varepsilon}

\newcommand{\R}{\mathbb R}

\newcommand{\Otr}{\Omega\times\R^k}

\renewcommand{\vec}[1]{\boldsymbol{#1}}

\begin{document}
\Large  \begin{center}{\bf Extended least action principle for
steady flows under a prescribed flux}\end{center} \normalsize
\begin{center} G. Wolansky \\ Department of Mathematics, \\
Technion, Haifa 32000, Israel \end{center}
\begin{abstract}
The extended principle of minimal action is described in the
presence of prescribed source and sink points. Under the
assumption of zero net flux, it leads to an optimal
Monge-Kantorovich transport problem of metric type. We concentrate
on action corresponding to a mecahnical Lagrangian. The optimal
solution turns out to be a measure supprted on a graph composed of
geodesic arcs connecting pairs of sources and sinks.
\end{abstract}
\section{Introduction}
\subsection{The extended action principle}
Consider the Lagrangian  of a mechanical system: $$L(p,x)=
\frac{1}{2}|p|^2 - V(x) \ , $$ where $x, p\in \R^k$ and $V$ is a
smooth bounded   potential function.

The {\it minimal} action is a function on $\R^k\times\R^k\times
\R\times\R$, defined by
\begin{equation}\label{classicJ}J(x_0,x_1,t_0,t_1):= \inf_{y}
\int_{t_0}^{t_1} L(\dot{y}(t), y(t)) dt \ ,  \end{equation} where
the infimum is taken on all orbits $y:[t_0,t_1]\rightarrow \R^k$
satisfying the end conditions $y(t_0)=x_0$, $y(t_1)=x_1$.
\par
The {\it extended minimal   action principle} (EMAP) was
introduced by Benamou and Brenier ([BB,BBG]), in the case $V\equiv
0$, and by [W] in the general case. Below we review its definition
and relation to the differential point of view of optimal mass
transportation - see also  ch. 8 of [V].
\par
Let $\overline{\cal M}$ be the set of probability Borel measures
on $\R^k$. Let $\vec{\cal E}$ the set of $\R^k-$valued Borel
measures defined on $\R^k$. The extended Lagrangian is defined for
any $(\mu, \overrightarrow{E})\in\overline{\cal M}\times \vec{\cal
E}$, provided $\overrightarrow{E}$ is absolutely continuous with
respect to $\mu$ and the Radon-Nykodym derivative
$d\overrightarrow{E}/d\mu\in \mathbb{L}^2(d\mu)$, as
$${\cal L}(\overrightarrow{E},\mu):= \frac{1}{2}\left|\frac{d \overrightarrow{E}}{d\mu}\right|^2-
V(x) \in \mathbb{L}^1(d\mu) . $$ The extended action is defined on
all orbits $\mu=\mu_t$ (res.
$\overrightarrow{E}=\overrightarrow{E}_t$) of $\overline{\cal M}$
(res. $ \overrightarrow{E}$) valued functions of the real line,
which satisfies the continuity equation
\begin{equation}\label{weakc} \partial_t\mu_t + \nabla_x\cdot
\overrightarrow{E}_t=0\end{equation} in the sense of
distributions.   For any pair $\mu_0,\mu_1\in \overline{\cal M}$
and any $t_0,t_1\in \R$ we define the EMAP
\begin{equation}\label{extedac}J(\mu_0, \mu_1, t_0,t_1):= \min_{(\mu,
\overrightarrow{E})}\int_{t_0}^{t_1}\int_{\R^k} {\cal L}\left(
\overrightarrow{E}_t, \mu_t\right)\mu_t(dx)dt\end{equation} over
all pairs $\left(\mu_t, \overrightarrow{E}_t\right)$ satisfying
(\ref{weakc}), subjected to the end condition $\mu_{t_0}=\mu_0$,
$\mu_{t_1}=\mu_1$.
\par
EMAP  is really an extension of the classical minimal action
principle, if the end measures $\mu_0,\mu_1$ are replaced by the
point measures $\delta_{(x_0)}$, $\delta_{(x_1)}$. In general, a
minimizer $(\overrightarrow{E}_t, \mu_t)$ of the EMAP exists and
 the Radon-Nykodym derivative $d \overrightarrow{E}_t/ d\mu_t$ is
a Lipschitz function on $Supp(\mu)$  under some general
conditions, so that the flow induced by
\begin{equation}\label{gradflow}\dot{x} = \left.\frac{d \overrightarrow{E}_t}{d\mu_t}
\right|_{x(t)}\end{equation} exists and is unique for $\mu$ a.e.
point $(x,t)$. Moreover, equation (\ref{gradflow}) is compatible
with the Euler-Lagrange equation  associated with the Lagrangian
$L$ for each individual orbit $x=x(t)$ (see [W]). \par The
motivation of Brenier and Benamou for the introduction of the
extended action principle is an attempt to devise an algorithm for
solving the Monge optimal mass transportation for quadratic cost:
Given two probability  measures $\mu_0,\mu_1$ on a common space
(say $\R^k$), find a mapping $T:\R^k\rightarrow \R^k$ which
transports  $\mu_0$ to $\mu_1$ so that the cost of transportation
 \begin{equation}\label{mk} \min_T\int_{\R^k} |T(x)-x|^2
d\mu_0 \ .  \end{equation}  Recall that a mapping $T$ transports
$\mu_0$ to $\mu_1$ ($T_\#\mu_0=\mu_1$) if and only if
$\mu_0\left(T^{-1}(A)\right)= \mu_1(A)$ for any (Borel) measurable
set $A$. The existence of a unique minimizer of (\ref{mk}) is
known if $\mu_0$ is absolutely continuous with respect to Lebesgue
measure and both $\mu_0,\mu_1$ have finite second moment (see [B],
and later extension in [GM]). It was the fundamental observation
of Benamou and  Brenier that (in the case $V\equiv 0$) the flow
(\ref{gradflow}) associated with the solution of the EMAP induces
the optimal mapping $T$ of (\ref{mk}) under some regularity
assumptions. The extension of this result to the case $V\not\equiv
0$, introduced in [W], relates the flow (\ref{gradflow}) to the
optimal solution of the mass transport with respect to the cost
function $c=c(x,y)$,
\begin{equation}\label{mkgen}
\min_T\int_{\R^k} c(x, T(x)) d\mu_0 \end{equation} where
$c(x,y)=J(x,y,t_1,t_2)$  given by the classical action
(\ref{classicJ}).

\vskip .2in
\subsection{Sources and sinks}
In this paper we consider the extended action principle under a
set of sources and sinks.
   Let  $D_i\subset \R^k$, $i=1,\ldots n$ be a disjoint set of
compact domains whose boundaries  $S_i=\partial D_i$ are  smooth
surfaces. On each such surface we assign an integrable function
$\lambda_i: S_i\rightarrow\R$ so that
\begin{equation}\label{tfc} \sum_{i=1}^N \oint_{S_i} \lambda_i d
s=0 \ . \end{equation}
\par
The extended action principle under the presence of prescribed
fluxes $\lambda_i$ across the surfaces $S_i$ is defined as
(\ref{extedac}), where the set of pairs $\left( \mu,
\overrightarrow{E}\right)$ is defined on $\left(\R^k-\cup_1^n
D_i\right)\times [t_0,t_1]$ and (\ref{weakc}) replaced by
\begin{equation}\label{weakcC}
\partial_t\mu + \nabla_x\cdot \overrightarrow{E}=0 \ \ \ on \ \
\left(\R^k-\cup_1^n D_i\right)\times [t_0,t_1]
 \ , \ \ \
\overrightarrow{E}_t\cdot\vec{n}_i= \lambda_i \ \ ; \ \ \ \forall
(s,t)\in S_i\times [t_0,t_1] \ , \end{equation} where
$\vec{n}_i(s)$ are the outward normal to $S_i$ at point $s$.
\par
We shall concentrate on  stationary  extended minimal action:
 $\mu_t$ ( and,
correspondingly, $\overrightarrow{E}_t$), are independent of $t$.
\par
Under stationarity condition, the extended principle of minimal
action takes the form \begin{equation}\label{minac} \min_{(\mu,
\overrightarrow{E})} \left[ \frac{1}{2} \int_{\R^k- \cup D_i}
\left|\frac{ d \overrightarrow{E}}{d\mu}\right|^2 d\mu -
\int_{\R^k- \cup D_i} V d\mu\right]\end{equation} where $\mu\in
\overline{\cal M}$ is supported on $\R^k-\cup D_i$ and $
\overrightarrow{E}$ is a vector valued measure, absolutely
continuous with respect to $\mu$, which satisfies in addition
\begin{equation}\label{add} \nabla\cdot\overrightarrow{E}=0 \ \ \
on \ \ \R^k-\bigcup D_i\ \ \ ; \ \ \ \overrightarrow{E}\cdot
\vec{n}_i=\lambda_i \ \ on \ S_i, \ 1\leq i\leq n \ .
\end{equation}
 \par\noindent{\bf Points sources and sinks}:  Letting the domains
 $D_j$ shrink to points $x_j\in\R^k$, we replace the functions $\lambda_i$ on
 $S_i$ by constants $\lambda_i\in \R$ so that
 $$ \sum_1^n \lambda_i=0 \ . $$
     The admissibility condition (\ref{add}) for the stationary minimal action
(\ref{minac}) can now be casted into the {\it single} condition:
\begin{equation}\label{div0} \nabla\cdot\overrightarrow{E}= \sum_1^n
\lambda_i\delta_{x_i} \ \ \ on\ \ \R^k \ .\end{equation} The
problem is now reduced to finding the minimizer $\mu\in
\overline{\cal M}$ of
\begin{equation}\label{minacpoint}\underline{\cal J}:=  \min_{(\mu,
\overrightarrow{E})} \left[ \frac{1}{2} \int_{\R^k} \left|\frac{ d
\overrightarrow{E}}{d\mu}\right|^2 d\mu - \int_{\R^k} V
d\mu\right]\end{equation} where $ \overrightarrow{E}$ is
absolutely continuous with respect to $\mu$ and satisfies
(\ref{div0}).
\subsection{Main results}\label{MainR}
In this paper we concentrate on point sources and sinks.  The set
of fluxes $\lambda_1, \ldots \lambda_n$ is divided into two sets:
sources ($\lambda_i>0$) and sinks ($\lambda_i<0$).
 \par Let \begin{equation}\label{Ipm} I_+:= \left\{ i \ ; \lambda_i>0\right\} \ \ \ ; \ \
\ I_-:= \left\{ i \ ; \lambda_i<0\right\}\ . \end{equation}
  \par
The {\it metric}  Monge transport plant between a pair of
probability measures $\mu_0,\mu_1$ is the special case of
(\ref{mkgen}) where $c(x,y)$ is a metric, say $c(x,y)=|x-y|$. It
turns out that the solution of the stationary EMAP is related to
the metric Monge  transport where $\mu_0,\mu_1$ are replaced by
the point measures \begin{equation}\label{dismu}
|\lambda|^{-1}\sum_{i\in I_+} \lambda_i\delta_{(x_i)} \ \ ; \
-|\lambda|^{-1}\sum_{i\in I_-}
\lambda_i\delta_{(x_i)}\end{equation} where $|\lambda|:=
\sum_{i\in I_+}\lambda_i\equiv -\sum_{j\in I_-} \lambda_j$. \par
In general, however, there are no  mappings $T$ which transport a
point measure $\mu_0$ to another measure, so there is no sense to
define the Monge problem (\ref{mkgen}) for this case. However, the
Monge problem can be relaxed to an optimization problem on the set
of 2-point probability distributions ${\cal Q}={\cal Q}(dxdy)$.
This is the celebrated Kantorovich  relaxation of the Monge
problem [K].
\par
We shall first describe the special case $V\equiv 0$,
corresponding to the Euclidian metric $c(x,y)=|x-y|$. In the case
of point measures (\ref{dismu}),  the Kantorovich relaxation of
the minimal cost (\ref{mkgen}) takes the form of
 the $1-$Wasserstein metric (see, e.g. [R]). It
is defined by a $|I_+|\times |I_-|$ matrix $A^o$ which realizes
the minimum
\begin{equation}\label{W10} {\bf W}^{(1)}(
\overrightarrow{\lambda})= \sum_{i\in I_+}\sum_{j\in I_-}
A^o_{i,j}|x_i-x_j|:= \min_{A\in {\cal
Q}_{\overrightarrow{\lambda}}} \sum_{i\in I_+}\sum_{j\in I_-}
A_{i,j}|x_i-x_j| \ .   \end{equation} where
\begin{equation}\label{calQ}{\cal Q}_{\overrightarrow{\lambda}}:=
\left\{A_{i,j}\geq 0 \ \ ; \ \ \sum_{j\in I_-} A_{i,j}=\lambda_i \
\  if \ \ i\in I_+ \ \  ; \ \ \sum_{i\in I_+} A_{i,j}=\lambda_j \
\  if \ \  j\in I_- \right\} \ . \end{equation} The first result
states that the minimal action (\ref{minacpoint}) is given by
$$\underline{\cal J}= \frac{1}{2}\left[{\bf W}^{(1)}\left(
\overrightarrow{\lambda}\right)\right]^2 \ . $$
\par
Let $G_0$ be the bi-graph composed of the set of vertices $x_i$,
$i=1,\ldots n$, and edges $L(i,j)$ are  segments connecting $i\in
I_+$ to $j\in I_-$.

 The second result states that the action minimizer of (\ref{div0}, \ref{minacpoint}) is
supported in the edges of the bi-graph $G_0$ for which
$A^o_{i,j}>0$. On each such edge, $\mu$ is a uniform  measure and
$$\mu(L(i,j)) = |\lambda|^{-1}A^o_{i,j} \ . $$
\par\noindent
We now consider the case $V\not\equiv 0$. We assume that
\begin{equation}\label{Vass}\lim_{|x|\rightarrow\infty}V(x)=0
 \ . \end{equation} Let
\begin{equation}\label{Vbar}\overline{V}=\sup_{\R^n}V\ \ . \end{equation} For any $E>\overline{V}$ we
consider the Riemannian metric associated with the Maupertuis'
action principle (see [Ar]): \begin{equation}\label{vdotdot}
d\sigma_E = \sqrt{E-V} ds \ , \end{equation} here $ds$ is the
Euclidean metric.
\par
 \par
 The geodesic distance  associated with this metric is denoted by $D_E(\cdot,\cdot)$.
We recall that a geodesic arc connecting two point $x,y$ coincides
with  an  orbit of the mechanical system $$ \ddot{x}+
\nabla_xV(x(t))=0 \ \  x(0)=x_0, \ x(T)=x_1 \ , $$ corresponding
to the energy level $|\dot{x}|^2/2 + V(x)=E$. Here
$T=T_{x_0,x_1}(E)$ is the time of flight from $x_0$ to $x_1$
(which is, of course, a function of $E$ as well).
\par
   The bi-graph $G_E$ is
defined, analogously to $G_0$, as the collection of vertices
$x_i$, $1\leq i\leq n$, and all edges composed of geodesic arcs
(with respect to the metric $d\sigma_E$) connecting $x_i$, $i\in
I_+$, to $x_j$, $j\in I_-$.  The $1-$Wasserstein metric associated
with this distance is given, analogously to ${\bf W}^{(1)}$, as
\begin{equation}\label{W1E} {\bf W}^{(1)}_E( \overrightarrow{\lambda})= \sum_{i\in
I_+}\sum_{j\in I_-} A^E_{i,j} D_E(x_i,x_j):= \min_{A\in {\cal
Q}_{\overrightarrow{\lambda}}} \sum_{i\in I_+}\sum_{j\in I_-}
A_{i,j} D_E(x_i,x_j) \ .   \end{equation} Then the minimal action
(\ref{minacpoint}) is
\begin{equation}\label{underA}\underline{\cal J}=\max_{E\geq
\overline{V}}\left[\sqrt{2}{\bf W}^{(1)}_E\left(
\overrightarrow{\lambda}\right)-E\right] \ . \end{equation} There
exists a minimizer $\mu\in \overline{\cal M}$ realizing this
action which satisfies the following: \par\noindent {\bf Case a:}
\   $E_0>\overline{V}$ is the maximizer of the RHS of
(\ref{underA}). Then there exists an action minimizer supported on
the bi-graph $G_E$ so that $\mu\left(L_E(i,j)\right)>0$ only if
$A^E_{i,j}>0$. To wit: \begin{equation}\label{wit}
\mu\left(L_{E_o}(i,j)\right)= 2^{-1/2}A_{i,j}^{E_o} T_{i,j}(E_o)\
\ , \end{equation} where $T_{i,j}(E_o)= T_{x_i,x_j}(E_o)$ as
defined below (\ref{vdotdot}). \vskip .2in\noindent {\bf Case b:}
If $E_0=\overline{V}$ is the maximizer of (\ref{underA}) then the
following holds: Let $\mu_0$ be the measure supported on
$G_{\overline{V}}$ subjected to (\ref{wit}). Then there exists
$\beta\in [0,1)$ so that $$ 2^{-1/2}\sum_{i\in I_+}\sum_{j\in
I_-}A_{i,j}^{\overline{V}} T_{i,j}(\overline{V})=1-\beta\leq 1 \ .
$$
 Let $x_0$ be a maximizer
of $V$, that is, $\overline{V}=V(x_0)$. Then  $$\mu=\mu_0+
\beta\delta_{x_0}$$ is an action minimizer. \par In both cases (a)
{\it and} (b), the following claim is valid:  \vskip .2in\noindent
{\bf Time/Flux duality:}  {\it The expectation of the inverse flow
time, $\mathbb{E}_\mu\left(T^{-1}\right)$, is proportional to the
total in(out) flux $|\lambda|:= \sum_{i\in I_+} \lambda_i =
-\sum_{j\in I_-} \lambda_j$:  }
\begin{equation}\label{timeflux}\mathbb{E}_\mu(T^{-1}):=\sum_{i\in I_+}\sum_{j\in I_-}
\mu\left( L_{E_0}(i,j)\right) T^{-1}_{i,j}(E_o) =
\frac{|\lambda|}{\sqrt{2}} \ . \end{equation} \vskip.2in\noindent
\subsection{Outline}
In section~\ref{weakfor} we derive the weak formulation of the
stationary EMAP, which leads to a dual problem:
$$ \inf_{\mu\in\ooM}\sup_{\phi\in C^1} \left[-\int_{\R^k-\cup D_i}
\left[V+\frac{1}{2}|\nabla\phi|^2\right]\mu(dx)
+\sum_i\oint_{S_i}\phi(s)\lambda_i(s) ds \right] \ . $$ In
section~\ref{pointss} we concentrate in the case of point sources
and sinks, where $D_i$ shrink to points $x_i$. It contains some
definitions and a preliminary lemma on the Wasserstein metric and
its dual representation.
\par
Section~\ref{proofmain} Is the most technical part of this paper.
It contains a sequence of auxiliary lemmas, which are needed to
the proof of the main result, as described in section~\ref{MainR}.
 The proof
itself is given at the end of this section. It is given for the
general case $V\not\equiv 0$, since the case $V\equiv 0$ follows
easily form the general one. Finally, the short section~\ref{conc}
summarizes the results of this paper.
\section{Weak formulation}\label{weakfor}
\subsection{Notations} Given $n$ disjoint compact sets $D_i\subset \R^k$ and $S_i\equiv
\partial D_i$ smooth surfaces, set $$\Omega:= \R^k-\bigcup_1^n D_i
\ . $$
 Let  $\overline{\cal M}$ stands for the
set of probability Borel measures on $\Omega$. Let
$${\cal M}:= \left\{ \nu \ ; \ \nu(dxdp) \ \text{Borel Probability measure on} \ \Omega\times\R^k \ , \ \
\int_{\Omega\times\R^{k}}|p|^2\nu(dxdp) < \infty \ . \right\}$$
  Let $ \mathbb{P}:\Otr
\rightarrow \Omega$ be the projection  $ \mathbb{P}(x,p)=x$. For
$\nu\in{\cal M}$,  $ \mathbb{P}_\# \nu$ is the push-forward
 to the marginal measure $\mu=\mu(dx)\in
\overline{\cal M}$:
$$\int_{\Otr} \psi(x)\nu(dxdp)= \int_{\Omega} \psi(x)
\mathbb{P}_\#\nu(dx)  \ \ \ ; \ \ \ \forall \psi\in C(\Otr) \ . $$
 Let also
$$\overrightarrow{E}_\nu(dx):=\int_{p\in\R^k} \vec{p}\nu(dpdx) \ \ .   $$
For $\mu\in\overline{\cal M}$ let $${\cal M}_\mu:= \left\{ \nu\in
{\cal M} \ ; \ \mathbb{P}_\#\nu=\mu\right\} \\ . $$
$$\Lambda:= \left\{ \nu\in {\cal M} \ ; \ \ \int_\Omega
\overrightarrow{E}_\nu(dx)\cdot\nabla\phi(x) -
\sum_i\oint_{S_i}\phi(s)\lambda_i(s) dpds = 0 \ \ \ \forall \phi \
\in C^1(\overline{\R^k} ) \right\}\ . $$

$$\Lambda_\mu:= \Lambda\cap {\cal M}_\mu \ .  $$

 The Lagrangian  is now defined as a function on ${\cal M}$ via:
  \begin{equation}\label{TE}
L(\nu)=\frac{1}{2}\int_{\Otr} |p|^2\nu(dxdp)-\int_\Omega
V(x)\mathbb{P}_\#(\nu)(dx) \ .
\end{equation}
\subsection{Weak form of the minimal action}\label{wfma}
We now describe the weak form of the stationary action principle:
 $$ \min_{\nu\in\Lambda} L(\nu) \ . \eqno{\bf (P)} $$
 Equivalently, if
 $$ \overline{L}(\mu):= \inf_{\nu\in \Lambda_\mu} L(\nu)$$
 then ($P$) is equivalent to
$$  \inf_{\mu\in\overline{\cal
M}}\overline{L}(\mu) \ . \eqno{({\bf P}^*) }$$ \subsection{Dual
representation} Let now define
    ${\cal L}:{\cal M}\times
C^\infty(\overline{\Omega}) \rightarrow \R\cup\{\infty\}$ by $$
{\cal L}(\nu,\phi):= L(\nu) - \int_\Omega
\overrightarrow{E}_\nu(dx)\cdot\nabla\phi(x) +
\sum_i\oint_{S_i}\phi\lambda_i  \ . $$
\par Next, we use an appropriate version of the minmax principle
 to obtain the {\it dual formulation}:
\begin{proposition}\label{prop1} For any $\mu\in\overline{\cal M}$,   $$\overline{L}(\mu)=\sup_{\phi\in
C^1(\overline{\Omega})}\inf_{\nu\in {\cal M}_\mu}{\cal
L}(\nu,\phi) .
$$
\end{proposition}
\begin{proof}
First, note that
$$ \overline{L}(\mu)=\inf_{\nu\in {\cal M}_\mu}\sup_{\phi\in
C^1(\overline{\Omega})}{\cal L}(\nu,\phi) \ . $$ Indeed, by
definition, $\sup_{\phi\in C^1(\overline{\Omega})}{\cal
L}(\nu,\phi)=\infty$ if $\nu\not\in\Lambda$, while ${\cal
L}(\nu,\phi)= L(\nu)$ otherwise. Next, note that $L$ is an affine
function  on each of the domains ${\cal M}_\mu$ and
$C^1(\overline{\Omega})$, separately. As such, it is a convex
functional on ${\cal M}_\mu$ and concave on
$C^1(\overline{\Omega})$. In addition, ${\cal M}_\mu$ is a compace
set with respect to the weak topology (in which ${\cal L}$ is
continuous). The Minmax theorem, then,  can  be applied (see, e.g.
[Ro]), and the claim follows.
\end{proof}
 Next we evaluate $$ \inf_{\nu\in {\cal
M}_\mu}{\cal L}(\nu,\phi)= \inf_{\nu\in {\cal
M}_\mu}\int_{\Otr}\left[ \frac{1}{2}|\vec{p}|^2 -
\vec{p}\cdot\nabla\phi -V\right]\nu(dxdp)
+\sum_i\oint_{S_i}\phi(s)\lambda_i(s) dpds  \ . $$ $$=
\inf_{\nu\in {\cal M}_\mu}\int_{\Otr}\left[
\frac{1}{2}|\vec{p}-\nabla\phi|^2 - \frac{1}{2}|\nabla
\phi|^2-V\right]\nu(dxdp) +\sum_i\oint_{S_i}\phi(s)\lambda_i(s)
dpds \ . $$ \begin{equation}\label{deltap} = -\int_{\R^k}\left[
\frac{1}{2}|\nabla\phi|^2+V\right]\mu(dx)
+\sum_i\oint_{S_i}\phi(s)\lambda_i(s) dpds  \ . \end{equation}
where the infimum is obtained at
$\nu(dxdp)=\mu(dx)\delta_{(\vec{p}-\nabla\phi)}$.
\par
Let us now define, for any $\mu\in \overline{\Omega}$ and $\phi\in
C^1(\overline{\Omega})$:
$$ J_\mu(\phi):= -\int_{\Omega}
\left[V+\frac{1}{2}|\nabla\phi|^2\right]\mu(dx)
+\sum_i\oint_{S_i}\phi(s)\lambda_i(s) ds \ . $$ By
Proposition~\ref{prop1} we get
\begin{corollary}\label{Ducor} \ \ \
$ \overline{L}(\mu)=\sup_{\phi\in
C^1(\overline{\Omega})}J_\mu(\phi)$
\end{corollary}
 \begin{proposition} If \ $\overline{\Omega}$ is compact, then there
 exists $\mu\in\overline{\cal M}$ which
 solves problem ${\bf
P}^*$.
\end{proposition}
\begin{proof}
If $\overline{\Omega}$ is compact, so is the set $\overline{\cal
M}$ with respect to the weak topology, as the set of Probability
Borel measures on a compact set. In addition, $\overline{L}$ is
lower semi continuous, since it is a supremum of the affine
functionals $J_{(\cdot)}(\phi)$ by Corollary~\ref{Ducor}. Hence, a
minimizing sequence of $\overline{L}$ in $\overline{\cal M}$
contains a subsequence which converges to a minimum of
$\overline{L}$.
\end{proof}
\section{Point sources and sinks}\label{pointss}
Assume now that the surfaces $S_i$ degenerate to points $x_i\in
\R^k$. The fluxes functions $\lambda_i$ defined on $S_i$
degenerate, then, to constants $\lambda_i\in \R$. The total flux
condition takes the form
\begin{equation}\label{tfcp} \sum_1^n \lambda_i=0 \ . \end{equation}
In this case, the functional $J_\mu$ is written as $$
J_\mu(\phi):= -\int_{\R^k}
\left[V+\frac{1}{2}|\nabla\phi|^2\right]\mu(dx)
+\sum_i\lambda_i\phi(x_i) \ . $$
\par\noindent

Recall that  $D_E$ be the distance  metric induced  by the
Riemannian  metric  $d\sigma_E$ (\ref{vdotdot}), that is, for
$x,x^{'}\in \R^k$: \begin{equation}\label{de} D_E(x, x^{'}):=
\inf_q\int_0^S \sqrt{E-V(q(s))}ds\end{equation} where the infimum
above is taken over all orbits $q(s)$ in Euclidian  arc-length
parameterization connecting $x$ to $x^{'}$. A  minimizer  orbit
$q$ of (\ref{de}) is called an $E-$ geodesic for the metric
$D_E(x, x^{'})$.
 \par
 Given $E>\overline{V}$ (\ref{Vbar}), $i,j\in (1, \ldots n)$, let
$q_{i,j}(s)$ be the $E-$geodesic  of  $D_E(x_i,x_j)$.
parameterized by the Euclidian arc-length. So, $q_{i,j}(0)=x_i$,
$q_{i,j}(S_{i,j})=x_j$ where $S_{i,j}$ is the Euclidian length of
the $D_E(x_i,x_j)$ geodesic curve. Let
\begin{equation}\label{Tijdef}T_{i,j}(E)=\int_0^{S_{i,j}} \left(E-V(q_{i,j}(s))\right)^{-1/2}
ds\ \ . \end{equation}  Note that $T_{i,j}(E)$ is the time
interval of existence of the orbit which connects $x_i$ to $x_j$
at energy $E$.  By differentiation of $D_E(x_i,x_j)$ with respect
to $E$ we obtain
\begin{lemma}\label{tijdef1} If $E>\overline{V}$ then
$$ \frac{d}{dE} D_E(x_i,x_j) = \frac{1}{2}T_{i,j}(E)  < \infty \ .
$$
\end{lemma}
 For each $E>\overline{V}$ define
\begin{equation}\label{rhoijdef}\rho_{i,j}(s):= \frac{1}{T_{i,j}(E)\sqrt{E-V(q_{i,j}(s))}}  \   ; 0\leq
s\leq S_{i,j} \ . \end{equation} Then $\rho_{i,j}$ is a
probability density on the interval $[0. S_{i,j}]$. Then
\begin{equation}\label{muijdef} \mu_{i,j}:= q_{i,j, \#}
\left(\rho_{i,j}(s)ds\right)\end{equation} is a push forward of
this probability density to a probability measure $\mu_{i,j}$
supported on the $E-$geodesic arc connecting $x_i$ and $x_j$ in
$\R^k$. It is defined, for any test function $\phi\in C(\R^k)$,
via $$ \int_{\R^k} \phi(x)d\mu_{i,j} = \int_0^{S_{i,j}}
\phi\left(q_{i,j}(s)\right) \rho_{i,j}(s)ds \ \ . $$

We now recall the definition of the  $1-$Wasserstein metric
(\ref{W1E}).
 By duality formulation of the
Wasserstein metric we also obtain (see [R]):
\begin{equation}\label{W1dual} {\bf W}_E^{(1)}(\lambda):= \max_{\overrightarrow{\phi}\in \R^n} \left\{
\sum_1^n\lambda_i\phi_i \ \ ; \ \ \frac{|\phi_j -\phi_k|}{D_E(x_j,
x_k)}  \leq 1 \ \ , \ \ \forall \ j,k\ \in \{1, \ldots n\}\ .
\right\} \ . \end{equation} Below we collect some useful results:
\begin{lemma}\label{du}
$A\in {\cal Q}_{ \overrightarrow{\lambda}}$ (\ref{calQ}) is an
optimal solution of (\ref{W1E}) iff there exists a minimizer
$\overrightarrow{\phi}$ of (\ref{W1dual}) so that  $A_{i,j}=0$ for
any pair $i,j$ satisfying $|\phi_i-\phi_j|< D_E(x_i,x_j)$.
\end{lemma}
\begin{proof}
For any $\overrightarrow{\phi}$ satisfying $\max_{i\not= j}
|\phi_i-\phi_j|/ D_E(x_i,x_j)\leq 1$ and $A\in{\cal Q}_{
\overrightarrow{\lambda}}$ we obtain
$$ \sum_{i\in I_+}\sum_{j\in I_+} A_{i,j} D_E(x_i,x_j) \geq
\sum_{i\in I_+}\sum_{j\in I+-} A_{i,j}(\phi_i-\phi_j) $$
$$= \sum_{i\in I_+}\left(\sum_{j\in I+-} A_{i,j} \phi_i\right)
-\sum_{j\in I_-}\left(\sum_{i\in I_+} A_{i,j} \phi_j\right) =
\sum_{i\in I_+} \lambda_i\phi_i - \sum_{j\in I_-}\lambda_j\phi_j =
\sum_{i=1}^n \lambda_i\phi_i \ , $$ where we used $A_{i,j}\equiv
0$ if either $i,j\in I_+$ or $i,j\in I_-$.  Since equality holds
for the optimal $A$ and $\overrightarrow{\phi}$, it follows that
$A_{i,j}>0$ implies $\phi_i-\phi_j= D_E(x_i,x_j)$.
\end{proof}
%\begin{lemma}\label{du1}
%\end{lemma}
 Recall the
definition of the subgradient of a function $h$ at $x\in \R^k$:
$$ \partial_x h= \left\{ p\in \R^k \ ; \ h(x^{'}) - h(x) \geq
p\cdot(x^{'}-x) \ \ \forall x^{'}\in\R^k \ . \right\}$$  The
following result can be found in, e.g., [HL]:\newpage
\begin{lemma}\label{ekland} \par\noindent
\begin{description}
\item{i)} \ If $h$ is a convex function defined on $\R^k$, then $\partial_xh\not=\emptyset$ for all $x\in \R^k$.
\item{ii)} \ If $h$ is convex, then $p\in \partial_x h$ if and
only if $x\in \partial_p h^*$ where $h^*$ is the Legendre
transform of $h$.
\item{iii)} \ If $h$ is convex, then for any $x,p\in \R^k$, the
inequality $h(x)+h^*(p)\geq x\cdot p$ holds with equality if and
only if $p\in\partial_xh$ if and only if $x\in \partial_p h^*$.
\end{description}
\end{lemma}

\section{Proof of the main result}\label{proofmain}
We shall  prove the main result of section~\ref{MainR}  for the
general case $V\not\equiv 0$ satisfying (\ref{Vass}). The special
case $V\equiv 0$   follows easily.
\par
 Let
$\overrightarrow{\phi}:= ( \phi_1, \ldots \phi_n)\in \R^n$ and
$$B_{\overrightarrow{\phi}}:=\left\{ \phi\in C^1(\R^k)
\ ; \ \phi(x_i)=\phi_i\right\} \ . $$ For $\mu\in \ooM$ set $$
H(\overrightarrow{\phi}; \mu):= \inf_{\phi\in
B_{\overrightarrow{\phi}}}\left[ \frac{1}{2}\int_{\R^k}
|\nabla\phi|^2 d\mu + \int_{\R^k} V d\mu\right] \ , $$ and $$
\overline{H}(\overrightarrow{\phi}):= \sup_{\mu\in \ooM} H(
\overrightarrow{\phi}; \mu) \  \ .  $$
\par\noindent
{\bf Remark:} \ Note that $H(\cdot; \mu)$ is not necessarily a
convex function on $\R^k$ for each $\mu\in\ooM$. However,
$\overline{H}$ is convex, as we shall see later on in
Corollary~\ref{corfinal}.
\begin{lemma}\label{6.3}
 For any $\overrightarrow{\phi}\in \R^n$,
\begin{equation}\label{dualvar} \overline{H}( \overrightarrow{\phi})
=
\min_{E \geq \overline{V}}\left \{
 E, \max_{i\not= j}\frac{|\phi_i-\phi_j|}{D_E(x_i,x_j)}\leq
 \sqrt{2}\right\}
\ .
\end{equation}
 \end{lemma}
\begin{proof}
 Let $E>\overline{V}$. Assume there exists $\phi\in C^1(\R^k)$ so that
 \begin{equation}\label{fff} |\nabla\phi|\leq \sqrt{2}\sqrt{E-V} \ \  on  \  \
 \overline{\R^k}\ and
 \ \ \phi(x_i)=\phi_i \ \ , \ \  1\leq i\leq n \ . \end{equation}
 Now, if such a function satisfies (\ref{fff}), then for any
 $i\not=j$ and any geodesic arc $y(t)$, $0\leq t\leq 1$,  connecting $x_i$ to $x_j$
 we obtain:
 $$ |\phi(x_i)-\phi(x_j)| \leq \int_0^1 \left|\nabla_{x(t)}\phi\right| \  \cdot \left|\dot{x}_{(t)}\right|dt \leq
 \sqrt{2}D_E(x_i,x_j) \ . $$
 In addition,  for any $\mu\in{\cal M}$
 $$ H(\overrightarrow{\phi}; \mu) \leq \frac{1}{2}\int_\R^k |\nabla\phi|^2
 d\mu + \int_\R^k V d\mu\leq E \ . $$
We now show the existence of $\phi$ satisfying (\ref{fff}),
provided \begin{equation}\label{maxi}\max_{i\not=
j}\frac{|\phi_i-\phi_j|}{D_E(x_i,x_j)} <
 \sqrt{2}\end{equation} Given $\eps>0$, let
 $D^\eps_E(x,x_i):= \max\left\{ D_E(x,x_i), \eps\right\}$.
 Set
 $$ \phi^\eps(x)=\min_{1\leq i\leq n} \left\{ \sqrt{2}D^\eps_E(x,x_i) +
 \phi_i-\eps\right\}  \ . $$
 Note that  $D^\eps_E(x, x_i)$ is a Lipschitz  function on $\R^k$  (recall $E>\overline{V}$) and satisfies
$|\nabla D_E(\cdot, x_i)|\leq \sqrt{E-V}$ for almost any
$x\in\R^k$.  In addition, the condition (\ref{maxi}) implies that
$\phi^\eps(x)= \phi_j$ if $|x-x_j|<\eps$ for any $j\in \{1, \ldots
n\}$, if $\eps$ is sufficiently small. Now, let us take the
smoothing kernel \begin{equation}\label{smk} \eta\in
C^\infty(\R^k) \ ; \ \ \eta\geq 0, \ \ ; \ \ \int_{\R^k}\eta = 1 \
, \ \ \ ; \ \ \eta(x)=0 \ if\ \ |x|>1 \ . \end{equation} Let
$\eta_\delta(x)= \delta^{-k}\eta(x/\delta)$. Let
$$\phi(x):= \eta_\delta * \phi^\eps \ . $$ Evidently, the
inequality $|\nabla\phi^\eps|\leq\sqrt{2}\sqrt{E-V}$ holds now
{\it everywhere} for $\phi$.  If $\delta<\eps$ then, by the last
condition in (\ref{smk}),  also the condition
$\phi^\eps(x_i)=\phi_i$ are preserved for $\phi$. So
\begin{equation}\label{firstineq}\overline{H}(\overrightarrow{\phi}) \leq  \inf_{E>
\overline{V}}\left \{
 E, \max_{i\not= j}\frac{|\phi_i-\phi_j|}{D_E(x_i,x_j)}\leq
 \sqrt{2}\right\}  \end{equation}
 is verified.
\par
To prove the opposite inequality we construct  a probability
measure as follows. Let $\overline{E}$ defined by
\begin{equation}\label{ebardef}\sqrt{2}D_{\overline{E}}(x_m, x_k)=|\phi_k-\phi_m| \ \ \ ; \
\ \ \sqrt{2}D_{\overline{E}}(x_i, x_j)\geq |\phi_i-\phi_j|
 \  \forall i\not= j \ . \end{equation}
  Assume $\overline{E} > \overline{V}$. Let $\mu=\mu_{k,m}$ as defined in (\ref{muijdef}). We now calculate the
minimizer $\phi\in B_{\overrightarrow{\phi}}$ of $\int_{\R^k}
|\nabla\phi|^2 d\mu_{km}$. Note that for any $\phi\in C^1(\R^k)$
we may define $\eta\in C^1\left[0,S_{k,m}\right]$ via
$\eta(s)=\phi\left(q_{k,m}(s)\right)$, where $q_{k,m}$ is the
parameterization of the geodesic arc connecting $x_m$ to $x_k$
(see paragraph preceding (\ref{Tijdef})). Then
\begin{equation}\label{phitoeta} \int_{\R^k} |\nabla\phi|^2 d\mu_{k,m} \geq
\int_0^{S_{k,m}} |\dot{\eta}(s)|^2 \rho_{k,m}(s) ds \
\end{equation}
where $\rho_{k,m}$ as defined in (\ref{rhoijdef}).  Now, the
minimizer on the RHS of (\ref{phitoeta}) subject to the condition
$\eta(0)=\phi_m$, $\eta\left(S_{k,m}\right)=\phi_k$  is attained
for
$$ \eta(s)= \phi_m + \frac{\phi_k-\phi_m}{\int_0^S\rho_{k,m}^{-1}}\int_0^s \frac{dt}{\rho_{k,m}(t)}  \ .  $$
By (\ref{phitoeta})  it follows that
\begin{equation}\label{firststep}\frac{1}{2}\int_{\R^k}
|\nabla\phi|^2 d\mu_{k,m}
\geq\frac{1}{2}\frac{|\phi_k-\phi_m|^2}{\int_0^{S_{k,m}}\rho_{k,m}^{-1}}
\ .
\end{equation} must holds for any $\phi\in B_{ \overrightarrow{\phi}}$.
\par
Next, by (\ref{ebardef}) and (\ref{rhoijdef}):
\begin{equation}\label{finstep}
\frac{|\phi_m-\phi_k|^2}{2\rho_{k,m}^2\left(\int_0^{S_{k,m}}\rho_{k,m}^{-1}\right)^2}
 + V = \overline{E} \ . \end{equation}
Multiply (\ref{finstep}) by $\rho_{k,m}$, integrate from $0$ to
$S_{k,m}$ and use (\ref{firststep}) to obtain $$ \int_\R^k \left[
\frac{|\nabla\phi|^2}{2}+ V\right] d\mu_{k,m} \geq \overline{E}$$
for any $\phi\in B_{ \overrightarrow{\phi}}$. This implies the
reverse inequality of (\ref{firstineq}), provided $\overline{E} >
\overline{V}$.
\par
Finally,  we observe that the choice $\mu= \delta_{x_0}$ where
$V(x_0)=\overline{V}$ guarantees: $$ \overline{H}(
\overrightarrow{\phi})\geq H(\overrightarrow{\phi}; \mu) =
\overline{V} \ . $$
 \end{proof}

Let now $\overline{H}^*$ be the Legendre transform of
$\overline{H}$:
 $$\overline{H}^*(\overrightarrow{\lambda}):= \sup_{
\overrightarrow{\phi}\in \R^n}\left\{ -\overline{H}(
\overrightarrow{\phi})  +\sum_1^n\lambda_i\phi_i\right\}
 \in \R\cup\{\infty\} \  . $$
Similarly
 $$H^*(\overrightarrow{\lambda};\mu):= \sup_{
\overrightarrow{\phi}\in \R^n}\left\{ -H(
\overrightarrow{\phi};\mu) +\sum_1^n\lambda_i\phi_i\right\}
  \in \R\cup\{\infty\} \ . $$ Note that {\it both} $\overline{H}^*$ and
 $H^*(\cdot;\mu)$ are convex by definition. From Lemma~\ref{6.3}
 we also obtain
 \begin{corollary}\label{meutar}
  $\overline{H}^*( \overrightarrow{\lambda})\in \R$ \  for each
 $\overrightarrow{\lambda}\in\R^k$ provided $\sum_1^n\lambda_i=0$.
 \end{corollary}
 \begin{proof}
We first note that $D_E(x_i,x_j)\approx \sqrt{E}$ as
$E\rightarrow\infty$. From Lemma~\ref{6.3} it follows that
$\overline{H}$ is at least quadratic with respect to
$|\phi_i-\phi_j|$ for any $i\not= j$. So, if we fix, say,
$\phi_1=0$, then
 $$\lim_{
 |\overrightarrow{\phi}|\rightarrow\infty, \phi_1=0}
 \frac{\overline{H}(\overrightarrow{\phi})}{|\overrightarrow{\phi}|}
 = \infty \ . $$
 In addition, $\overline{H}$ is invariant under the shift $\overrightarrow{\phi}\rightarrow \overrightarrow{\phi}+ t
 \overrightarrow{1}$  where $ \overrightarrow{1}=
 (1, \ldots 1)\in \R^k$ and $t\in \R$.
 This implies that $\sum_1^n\lambda_i\phi_i-\overline{H}(
 \overrightarrow{\phi})$ is bounded from
 above on $\R^k$, provided $\sum_1^n\lambda_i=0$.
 \end{proof}
From Proposition~\ref{prop1} and (\ref{deltap}) we observe that
\begin{equation}\label{h*} \overline{L}(\mu)= H^*( \overrightarrow{\lambda}; \mu) \ .
\end{equation}
By definition, $H^*( \overrightarrow{\lambda};\mu) \geq
\overline{H}^*( \overrightarrow{\lambda})$ for any $\mu\in \ooM$.
So, we obtain: \begin{corollary}\label{corh*} The minimal action
is not smaller than $\overline{H}^*( \overrightarrow{\lambda})$.
\end{corollary}
We now state
\begin{lemma}\label{lem6.2}
Let  $\overrightarrow{\phi}\in
\partial_{\overrightarrow{\lambda}}H^*(\cdot;\mu)$. Assume
$H(\overrightarrow{\phi};
\mu)=\overline{H}(\overrightarrow{\phi})$ and $H(\cdot;\mu)$ is a
convex function in the first variable. Then $\mu$ is an action
minimizer (that is, a minimizer of (\ref{h*})).
\end{lemma}
\begin{proof}
   Consider the
chain of inequalities:
$$ H^*(\overrightarrow{\lambda};\mu)+
\overline{H}(\overrightarrow{\phi}) \geq
H^*(\overrightarrow{\lambda};\mu) + H(\overrightarrow{\phi};\mu)
\geq \overrightarrow{\phi}\cdot \overrightarrow{\lambda} \ , $$
which holds for any $\mu\in\overline{\cal M}$,
$\overrightarrow{\phi}\in \R^n$ and  $\overrightarrow{\lambda}$ in
the domain of $\overline{H}^*$. The left inequality follows from
the definition of $\overline{H}$ as a maximizer of $H$ over
$\overline{\cal M}$, while the right inequality follows by the
convexity of $H(\cdot;\mu)$ and Lemma~\ref{ekland}-iii. The
conditions of the Lemma and Lemma~\ref{ekland}-iii again imply
that both inequalities are, in fact, equalities , so
$$ H^*(\overrightarrow{\lambda};\mu)+
\overline{H}(\overrightarrow{\phi}) =
 \overrightarrow{\phi}\cdot \overrightarrow{\lambda} \ . $$
 Now, by definition \begin{equation}\label{111}H^*(\overrightarrow{\lambda};\mu)\geq
 \overline{H}^*(\overrightarrow{\lambda})  \end{equation}
 so
 \begin{equation}\label{222} \overline{H}^*(\overrightarrow{\lambda})+
\overline{H}(\overrightarrow{\phi}) \leq
 \overrightarrow{\phi}\cdot \overrightarrow{\lambda} \ .
 \end{equation}
 Now, $\overline{H}^{**}\leq \overline{H}$ by definition, so
\begin{equation}\label{333} \overline{H}^*(\overrightarrow{\lambda})+
\overline{H}^{**}(\overrightarrow{\phi}) \leq
 \overrightarrow{\phi}\cdot \overrightarrow{\lambda}
 \end{equation}
holds as well.
 Now, Lemma~\ref{ekland}-iii implies that the {\it reverse} inequality must hold in (\ref{333}). Hence
there must be an equality in (\ref{333}), which induces the
equalities in (\ref{222}) and (\ref{111}) as well. It follows
  that $\mu$ is an action minimizer by (\ref{h*}) and
  Corollary~\ref{corh*}.
\end{proof}
 \begin{corollary}\label{cormain} \ \ \ \
$ \overline{H}^*( \overrightarrow{\lambda})
 =
 \sup_{E\geq \overline{V}}\left\{\sqrt{2}{\bf
 W}^{(1)}_E(\overrightarrow{\lambda})-E\right\} $ .
 \end{corollary}
 \begin{proof}
 By definition of $\overline{H}^*$ and Lemma~\ref{6.3},
 $$ \overline{H}^*(\overrightarrow{\lambda})=
 -\min_{\overrightarrow{\phi}}\left[
 \overline{H}(\overrightarrow{\phi})-
 \overrightarrow{\lambda}\cdot \overrightarrow{\phi}\right] =
-\min_{\overrightarrow{\phi}}\inf_{E\geq \overline{V}}\left[E -
\overrightarrow{\lambda}\cdot \overrightarrow{\phi} \ ; \
\max_{i\not=j}\frac{|\phi_i-\phi_j|}{D_E(x_i,x_j)}\leq
\sqrt{2}\right]
 $$
 \begin{equation}\label{W1e} = \sup_{E\geq \overline{V}}\left[ \sqrt{2}{\bf
 W}^{(1)}_E(\overrightarrow{\lambda})-E\right] , \end{equation}
 where we used the duality relation given by (\ref{W1dual}).
\end{proof}
 \begin{lemma}\label{mongk}
 Suppose $E_0>\overline{V}$ is the minimizer of  (\ref{W1e}). Then there exists $A^{E_0}\in {\cal Q}_{\overrightarrow{\lambda}}$ which minimize
 the Wasserstein  cost ${\bf W}^{(1)}_{E_0}$, i.e.  $\sum_{i\in I_+}\sum_{j\in I_-} A^{E_0}_{i,j}D_{E_0}(x_i,x_j)=
 {\bf W}^{(1)}_{E_0}( \overrightarrow{\lambda})$,   and
$$ \sum_{i\in I_+}\sum_{j\in I_-}
A^{E_0}_{i,j}T_{i,j}(E_0)=\sqrt{2} \  $$  is satisfied. If
$E_0=\overline{V}$ then for any such $A^{\overline{V}}\in {\cal
Q}_{\overrightarrow{\lambda}}$, the inequality $$ \sum_{i\in
I_+}\sum_{j\in I_-}
A^{\overline{V}}_{i,j}T_{i,j}(\overline{V})\leq \sqrt{2} \  $$
holds.
 \end{lemma}
\begin{proof}
Let $E_n\searrow E_0$. For each $n$, set $A^{E_n}\in {\cal
Q}_{\overrightarrow{\lambda}}$ be a minimizer of the
Monge-Kantorovich problem associated with ${\bf
W}^{(1)}_{E_n}(\overrightarrow{\lambda})$. Note that such a
minimizer may not be unique. We choose a subsequence so that the
limit
\begin{equation}\label{limae} A^{E_0^+}:= \lim_{n\rightarrow\infty} A_n^{E_n}\end{equation} exists.
Evidently,  $A^{E_0^+}$ is a minimizer of the Monge-Kantorovich
problem associated with ${\bf
W}^{(1)}_{E_0}(\overrightarrow{\lambda})$ (again, possibly not the
only one).
\par
Next, since $E_0$ is a maximizer of (\ref{W1e}),
$$ \sqrt{2}\left({\bf
 W}^{(1)}_{E_n}(\overrightarrow{\lambda})-{\bf
 W}^{(1)}_{E_0}(\overrightarrow{\lambda})\right)\leq E_n-E_0 \ , $$
 so
 $$ \sum_{i\in I_+}\sum_{j\in I_-}A^{E_n}_{i,j} \left(
 D_{E_n}(x_i,x_j) - D_{E_0}(x_i,x_j)\right)
 \leq \sum_{i\in I_+}\sum_{j\in I_-}A^{E_n}_{i,j}
 D_{E_n}(x_i,x_j) - \sum_{i\in I_+}\sum_{j\in
 I_-}A^{E^+_0}_{i,j}D_{E_0}(x_i,x_j)$$
 \begin{equation}\label{lastline} = {\bf
 W}^{(1)}_{E_n}(\overrightarrow{\lambda})-{\bf
 W}^{(1)}_{E_0}(\overrightarrow{\lambda})\leq \frac{E_n-E_0}{\sqrt{2}} \ .
 \end{equation}
 Take the limit $n\rightarrow\infty$ and use (\ref{limae}) and Lemma~\ref{tijdef1} to obtain from
 (\ref{lastline})
 \begin{equation}\label{show1} \sum_{i\in I_+}\sum_{j\in I_-}A^{E_0^+}_{i,j} T_{i,j}(E_0)\leq \sqrt{2} \ .
\end{equation}
This completes the proof for the case $E_0=\overline{V}$. \par
Now, let $E_n\nearrow E_0$, and let us consider the subsequence
along which the limit
\begin{equation}\label{limae1} A^{E_0^-}:= \lim_{n\rightarrow\infty} A_n^{E_n}\end{equation}
exists. By following the preceding argument we obtain
\begin{equation}\label{show2}  \sum_{i\in I_+}\sum_{j\in
I_-}A^{E_0^-}_{i,j} T_{i,j}(E_0)\geq \sqrt{2} \ .
\end{equation}
Finally, if the maximizer $A^{E_0}$ of ${\bf
W}^{(1)}(\overrightarrow{\lambda})$ is unique, then
$A^{E_0^-}_{i,j}=A^{E_0^+}_{i,j}:=A^{E_0}$ and the proof follows
from (\ref{show1}, \ref{show2}).  Otherwise, since both
$A^{E_0^\pm}_{i,j}$ are minimizers, so is the convex combination
thereof. Now, we utilize (\ref{show1}, \ref{show2}) to choose
$\alpha\in[0,1]$ for which $A^{E_0}:=\alpha A^{E_0^-}_{i,j} +
(1-\alpha)A^{E_0^+}_{i,j}$ satisfies the desired equality
$$  {\bf W}_{E_0}^{(1)}(\overrightarrow{\lambda}) \equiv\sum_{i\in I_+}\sum_{j\in
I_-}A^{E_0}_{i,j} T_{i,j}(E_0)= \sqrt{2} \ . $$
\end{proof}
Given $\overrightarrow{\phi}\in \R^n$ so that $\overline{H}(
\overrightarrow{\phi})=E$. Let \begin{equation}\label{IsubE}I_E:=
\left\{ (i,j)\in I_+\times I_- \ \ ; \ \ |\phi_i-\phi_j|=
\sqrt{2}D_E(x_i,x_j)\right\} \ .
\end{equation}
Let $C_E$ be the convex hall of the set of measures $\mu_{i,j}$ as
defined in (\ref{muijdef}) where $(i,j)\in I_E$, that is
\begin{equation}\label{661}\mu\in C_E \ \ \text{iff} \ \mu=\sum_{i\in I_+}\sum_{j\in I_-} \alpha_{i,j} \mu_{i,j} \ \
\ \ ; \ \ \ \alpha_{i,j}\geq 0 \  \ \ , \ \ \sum_{i\in
I_+}\sum_{j\in I_-} \alpha_{i,j}=1  \ \ \ . \end{equation}
\begin{lemma}\label{6.4}
For any $\overrightarrow{\phi}\in \R^n$ satisfying
$\overline{H}(\overrightarrow{\phi})=E>\overline{V}$
$$ \overline{H}( \overrightarrow{\phi})= H\left(\overrightarrow{\phi};
\mu\right)=E$$ holds for any $\mu\in C_E$.
 \end{lemma}
\begin{proof}
The proof of Lemma~\ref{6.4} is, basically, identical to the proof
of the second part of Lemma~\ref{6.3}.
\end{proof}
\begin{lemma}\label{6.6}
Let $E>\overline{V}$ and  $\mu\in C_E$ for some
$\overrightarrow{\phi}\in \R^n$ satisfying $\overline{H}(
\overrightarrow{\phi})=H(\overrightarrow{\phi};\mu)=E$. Then for
any $\ \overrightarrow{\zeta}\in \R^k$,
$$ H(\overrightarrow{\zeta};\mu)= \frac{1}{2}\sum_{i\in
I_+}\sum_{j\in I_-}
\alpha_{i,j}\frac{(\zeta_i-\zeta_j)^2}{D_E(x_i,x_j) T_{i,j}(E)} +
\int V d\mu \ . $$ In particular, $H(\cdot,\mu)$ is convex for any
$\mu$ of the form (\ref{661}).
\end{lemma}
\begin{proof}
Given $\overrightarrow{\zeta}\in \R^n$, let $\zeta\in C^1(\R^k)$
be an optimal solution corresponding to $\mu$.  We push it
backward to a function on the graph composed of $\cup_{i,j} [0,
S_{i,j}]$ via
$$\eta_{i,j}(s):=\zeta\left(q_{i,j}(s)\right) , $$
so \begin{equation} H(\overrightarrow{\zeta}; \mu)\leq
\frac{1}{2}\sum_{i,j}\alpha_{i,j}\int_0^{S_{i,j}} \rho_{i,j}(s)
|\dot{\eta}_{i,j}|^2 ds  + \int V d\mu \ .
\label{ttt}\end{equation} The equality in (\ref{ttt}) is achieved
if  we minimize the RHS on each branch separately, subjected to
the prescribed end conditions $\eta_{i,j}(0)= \zeta_i$,
$\eta_{i,j}(S_{i,j})=\zeta_j$.  Hence $\rho_{i,j}\dot{\eta}_{i,j}$
is constant on each branch. Taking the end conditions and the
definition of $\rho_{i,j}$ (\ref{rhoijdef}) we obtain that
\begin{equation}\label{rhoetafixed}\rho_{i,j}(s)\dot{\eta}_{i,j}(s)=\left(\int_0^{S_{i,j}}
\frac{ds}{\rho_{i,j}}\right)^{-1}(\zeta_j-\zeta_i) =
\frac{\zeta_j-\zeta_i}{T_{i,j}D_E(x_i,x_j)} \ . \end{equation}
Now, we perturb $\zeta_i\rightarrow \zeta_i+\eps$. Let
$\hat{\zeta}_{i,j}$ be defined on the $(i,j)$ branch $[0,
S_{i,j}]$ where $\hat{\zeta}(0)=1$, $\hat{\zeta}(S_{i,j})=0$. Let
$e_i$ be the canonical $i-$ unit vector in $\R^n$, then,
evaluating $H( \overrightarrow{\zeta}+ \eps e_i ; \mu)$ via the
function $\eta_{i,j}+\eps\hat{\zeta}_{i,j}$ on each branch $[0,
S_{i,j}]$ where $\alpha_{i,j}>0$ (fixed $i$) and integration by
parts yields
$$ H( \overrightarrow{\zeta}+ \eps e_i ; \mu) - H( \overrightarrow{\zeta} ; \mu)
= \eps \sum_j\alpha_{i,j} \rho_{i,j}(0)\dot{\eta}_{i,j}(0) +
O(\eps^2) \ , $$ From (\ref{rhoetafixed}) it follows that
$$ \frac{\partial H(\
\overrightarrow{\zeta}; \mu)}{\partial \zeta_i}=\left\{
\begin{array}{cc}
 \sum_{j\in I_-}
\alpha_{i,j}\frac{\zeta_i-\zeta_j}{D_E(x_i,x_j)T_{i,j}}  & i\in
I_+
 \\ \sum_{j\in I_+}
\alpha_{i,j}\frac{\zeta_i-\zeta_j}{D_E(x_i,x_j)T_{i,j}}
   & i\in I_-
\end{array} \right. $$
and the proof follows by integration.
\end{proof}
\begin{corollary}\label{corfinal}
The function $\overline{H}$ is convex on $\R^k$.
\end{corollary}
\begin{proof}
By definition, $\overline{H}( \overrightarrow{\phi})$ is the
maximum of the set $H(\overrightarrow{\phi};\mu)$ where $\mu$ run
on $\ooM$.
 From Lemma~\ref{6.6} it
follows that this maximum is obtained  at a convex function.
Hence, $\overline{H}$ is the maximum of a family of convex
functions, so it is convex.
\end{proof} \vskip .2in\noindent  {\bf Proof of the Main Result:} \\
\par We prove the main result in full generality ($V\not\equiv 0$).
 \par\noindent
We obtain from (\ref{h*}) and Corollary~\ref{cormain} that
(\ref{underA})  is a lower bound for the minimal action.
\par
By Corollary~\ref{meutar}, $\overrightarrow{\lambda}$ is in the
domain of $\overline{H}^*$ if $\sum_1^n\lambda_i=0$, so there
exists   $\overrightarrow{\phi}
\in\partial_{\overrightarrow{\lambda}}\overline{H}^*$. Let
$\overline{H}( \overrightarrow{\phi})=E$ and assume that
$E>\overline{V}$. By Lemma~\ref{6.3},
\begin{equation}\label{res} \max_{i\not= j}
\frac{|\phi_i-\phi_j|}{D_E(x_i,x_j)} \leq \sqrt{2} \ .
\end{equation} Since $\overline{H}^*( \overrightarrow{\lambda})+E
=\overline{H}( \overrightarrow{\phi})+ \overline{H}^*(
\overrightarrow{\lambda})=\sum\lambda_i\phi_i$ we obtain the
equality in (\ref{doo}) below from  Corollary~\ref{cormain}:
\begin{equation}\label{doo}{\bf
W}^{(1)}_E(\overrightarrow{\lambda})=2^{-1/2} \sum_1^n
\lambda_i\phi_i  \ .
\end{equation}
\par
Let \begin{equation}\label{mu} \mu=2^{-1/2}\sum_{i\in
I_+}\sum_{j\in I_-} A^E_{i,j} T_{i,j}(E) \mu_{i,j}\end{equation}
where $\mu_{i,j}$ as defined in (\ref{rhoijdef}, \ref{muijdef}),
corresponding to the same energy $E$. By Lemma~\ref{mongk}, $\mu$
is a convex combination of $\mu_{i,j}$, hence a  probability
measure.
\par
 Now let $A\in{\cal
Q}_{\overrightarrow{\lambda}}$. By (\ref{res}, \ref{doo})
$$ \sum_{i\in I_+}\sum_{j\in I_-} A_{i,j}D_E(x_i,x_j) \geq
2^{-1/2} \sum_{i\in I_+}\sum_{j\in I_-} A_{i,j}|\phi_i-\phi_j| $$
$$
\geq 2^{-1/2} \sum_{i\in I_+}\sum_{j\in I_-}
A_{i,j}(\phi_i-\phi_j)= 2^{-1/2}\sum_1^n \lambda_i\phi_i = {\bf
W}^{(1)}_E(\overrightarrow{\lambda}) \ . $$ Now substitute the
minimizer $A^E\in {\cal Q}_{ \overrightarrow{\lambda}}$ for $A$
above. Then the chain of inequalities turn into equalities. In
particular we obtain that, for any $i\in I_+$, $j\in I_-$, either
$A^E_{i,j}=0$ or $\phi_i-\phi_j = \sqrt{2}D_E(x_i,x_j)$. As a
result we can apply Lemma~\ref{6.4} to obtain
\begin{equation}\label{stop1} H(\overrightarrow{\phi};\mu)=E \ .
\end{equation}
 In addition, we
substitute the equalities $\phi_i-\phi_j = \sqrt{2}D_E(x_i,x_j)$
and $\alpha_{i,j}=2^{-1/2} A_{i,j}^E T_{i,j}(E)$ in
Lemma~\ref{6.6} to obtain $$ \nabla_{ \overrightarrow{\phi}}
H(\cdot, \mu)=\overrightarrow{\lambda}\ ,
$$ which implies
\begin{equation}\label{stop2}\overrightarrow{\phi}\in\partial_{\overrightarrow{\lambda}}H^*\left(\cdot;\mu\right)
\ \end{equation} By (\ref{stop1},\ref{stop2}) and
Lemma~\ref{lem6.2} we obtain that $\mu$ is an action minimizer. In
particular, it follows that (\ref{underA}) {\it is} the minimal
action.
\par
Finally, if $E=\overline{V}$, then let $\beta=2^{-1/2}\sum_{i\in
I_+}\sum_{j\in I_-} A^E_{i,j} T_{i,j}(E) $. By lemma~\ref{mongk},
$\beta\leq 1$. define
\begin{equation}\label{mu1} \mu=2^{-1/2}\sum_{i\in
I_+}\sum_{j\in I_-} A^E_{i,j} T_{i,j}(E) \mu_{i,j} +
(1-\beta)\delta_{x_0}\end{equation} where $x_0$ is a maximizer of
$V$, so $V(x_0)=\overline{V}$. Since $E=\overline{V}$ the equality
(\ref{stop1}) holds for $\mu$ given by (\ref{mu1}). In addition,
(\ref{stop2}) is also verified for this $\mu$ by Lemma~\ref{6.6}.
Hence, $\mu$ is an action  minimizer via Lemma~\ref{lem6.2} as
well.
\section{Conclusion}\label{conc}
We considered the extended minimal action principle for stationary
actions in the presence of point sources and sinks. The main
conclusion of this paper is that this minimal, stationary action
is obtained as a minimization of a {\it metric} Monge-Kantorovich
for the (non-normalized) pair of discrete measures (\ref{dismu}).
This stands in contrast to the non-metric Monge-Kantorovich
transport (\ref{mk},\ref{mkgen}) which is related to the extended
minimal action in the non-stationary case.
\par
Another interesting conclusion is the relation (\ref{timeflux})
between the expectation of the inverse-time of the minimal orbits
to the given flux. This equality follows from substitution of
(\ref{wit}) and using the fact $A^E\in {\cal Q}_{
\overrightarrow{\lambda}}$ in the middle term of (\ref{timeflux}).
We stress that this relation is preserved in the case
$E=\overline{V}$, since the point measure $\beta\delta_{(x_0)}$
corresponds to orbits of infinite time length-these are the orbits
which converges to the maximum of $V$, but {\it never} get there,
since the metric $\sqrt{E-V}ds$ is degenerate at the point $x_0$
where $V(x_0)=\overline{V}$.

\begin{center}{\bf References}\end{center}
\begin{description}
\item{[Ar]} \ V.I.Arnold, {\it Mathematical Methods of Classical
Mechanics}, Springer-Verlag, 1980
\item{[B]} \ Y. Brenier: {\it Polar factorization and monotone
rearrangement of vector valued  functions}, Arch. Rational Mech
\&Anal., {\bf 122}, (1993), 323-351.
\item{[BB]} \ J.D.Benamou, Y. Brenier: {\it A computational fluid
mechanics solution to the Monge-Kantorovich mass transfer
problem}, Numer.Math., {\bf 84} (2000), 375-393.
\item{[BBG]} \ J.D.Benamou, Y. Brenier and K.Guitter: {\it The Monge-Kantorovich mass transfer
and its  computational fluid mechanics formulation}, Inter. J.
Numer.Meth.Fluids, {\bf 40} (2002), 21-30. \item{[W]} \ G.
Wolansky, {\it Optimal Transportation in the presence of a
prescribed pressure field}, preprint
\item{[GM]} \ W. Gangbo and\& R.J. McCann: \ {\it The geometry of
optimal transportation}, Acta Math., 177 (1996), 113-161
\item{[HL]} J.-B Hiriart-Urruty, C. Lemarechal
\textit{Convex Analysis and Minimization Algorithms
II},Grundlehren der Mathematischen Wissenschaften, vol. 306,
Springr-Verlag, Berlin, 1993
\item{[K]} \ L. Kantorovich: \ {\it On the translocation of
masses}, C.R (Doclady) Acad. Sci. URSS (N.S), 37, (1942), 199-201
\item{[R]} \ S.T Rachev and L. R$\ddot{u}$schendorf: \ {\it Mass
Transportation Problems}, Vol 1, Springer, 1998 \item{[Ro]} \ R.T.
Rockafeller, {\it Convex Analysis}, \ Princeton, N.J, Princeton U.
Press
\item{[V]}\ C. Villani \ {\it Topics in Optimal Transportation},
Graduate Studies in Mathematics, Vol. 58, AMS
\end{description}
\end{document}